\documentclass[12pt]{amsart}
\usepackage{epsfig}

\usepackage{pb-diagram}  

\def\C{\mathbf{C}}
\def\F{\mathbf{F}}

\def\Z{\mathbf{Z}}

\def\G{\mathbf{G}}    

\def\SL#1{\mathbf{SL}_{#1}}

\def\a{\alpha}
\def\b{\beta}

\def\l{\lambda}

\def\s{\sigma}

\def\rmi{\uppercase\expandafter{\romannumeral1}}
\def\rmii{\uppercase\expandafter{\romannumeral2}}
\def\rmiii{\uppercase\expandafter{\romannumeral3}}
\def\rmv{\uppercase\expandafter{\romannumeral5}}
\def\rmvi{\uppercase\expandafter{\romannumeral6}}
\def\rmvii{\uppercase\expandafter{\romannumeral7}}
\def\rmviii{\uppercase\expandafter{\romannumeral8}}

\def\FF{\mathcal{F}}          
\def\M {\mathcal{M}}             
\def\O{\mathcal{O}}              

\def\RR{\mathcal{R}}
\def\W{\mathcal{W}}              

\def\g{\mathfrak{g}}          
\def\n{\mathfrak{n}}

\def\so{\mathfrak{so}}
\def\s{\mathfrak{s}}
\def\h{\mathfrak{h}}          
\def\Liesl{\mathfrak{sl}}     
\def\gl{\mathfrak{gl}}        
\def\so{\mathfrak{so}}      

\def\can#1{\left\langle#1\right\rangle}

\def\PB{\left\{\cdot\,,\cdot\right\}}

\def\pb#1{\left\{#1\right\}}

\def\set#1{\left\{#1\right\}}

\def\lb#1{\[#1\]}

\def\inn#1#2{\left\langle#1\,\vert\,#2\right\rangle}

\def\INN{\langle\cdot\,\vert\,\cdot\rangle}

\def\({\left(}
\def\){\right)}
\def\[{\left[}
\def\]{\right]}

\def\ad{\mathop{\rm ad}\nolimits}
\def\Ad{\mathop{\rm Ad}\nolimits}

\def\det{\mathop{\rm det}\nolimits}

\def\diff{\mathsf{d}}

\def\Trace{\mathop{\rm Trace}\nolimits}

\newenvironment{eqn*}[1][1.5]
  {$$\renewcommand{\arraystretch}{#1}
      \begin{array}{rcl}}
      {\end{array}$$}

\newtheorem{thm}{Theorem}[section]
\newtheorem{proposition}[thm]{Proposition}
\newtheorem{corollary}[thm]{Corollary}
\newtheorem{lemma}[thm]{Lemma}
\theoremstyle{definition}
\newtheorem{definition}[thm]{Definition}

\newenvironment{matr}[1]{\left(\begin{array}{#1}}{\end{array}\right)}

\def\ds{\displaystyle}
\def\comment#1{}  

\def\mi{{-1}} 

\def\i{\imath}
\def\p{\partial}
\def\pp#1#2{\frac{\p #1}{\p #2}}

\def\we{\wedge}

\def\transp#1{{#1^\top}} 
\def\wght{\varpi} 
\def\defi#1{\index{#1}\emph{#1}} 

\DeclareMathAlphabet{\mat}{U}{msb}{m}{n}
\DeclareMathAlphabet{\cur}{U}{eus}{m}{n}
\DeclareMathAlphabet{\mathfrak}{U}{euf}{m}{n}

\newcommand{\lra}{\longrightarrow}

\catcode`\á=13 \defá{\'a}
\catcode`\à=13 \defà{\`a}
\catcode`\â=13 \defâ{\^a}
\catcode`\ä=13 \defä{\"a}
\catcode`\é=13 \defé{\'e}
\catcode`\è=13 \defè{\`e}
\catcode`\ê=13 \defê{\^e}
\catcode`\ë=13 \defë{\"e}
\catcode`\í=13 \defí{\'\i}
\catcode`\ì=13 \defì{\`\i}
\catcode`\î=13 \defî{\^\i}
\catcode`\ï=13 \defï{\"\i}
\catcode`\ó=13 \defó{\'o}
\catcode`\ò=13 \defò{\`o}
\catcode`\ô=13 \defô{\^o}
\catcode`\ö=13 \defö{\"o}
\catcode`\ú=13 \defú{\'u}
\catcode`\ù=13 \defù{\`u}
\catcode`\û=13 \defû{\^u}
\catcode`\ü=13 \defü{\"u}

\catcode`\ç=13\defç{\c c}

\catcode`\Á=13 \defÁ{\'A}
\catcode`\À=13 \defÀ{\`A}
\catcode`\Â=13 \defÂ{\^A}
\catcode`\Ä=13 \defÄ{\"A}
\catcode`\É=13 \defÉ{\'E}
\catcode`\È=13 \defÈ{\`E}
\catcode`\Ê=13 \defÊ{\^E}
\catcode`\Ë=13 \defË{\"E}
\catcode`\Í=13 \defÍ{\'I}
\catcode`\Ì=13 \defÌ{\`I}
\catcode`\Î=13 \defÎ{\^I}
\catcode`\Ï=13 \defÏ{\"I}
\catcode`\Ó=13 \defÓ{\'O}
\catcode`\Ò=13 \defÒ{\`O}
\catcode`\Ô=13 \defÔ{\^O}
\catcode`\Ö=13 \defÖ{\"O}
\catcode`\Ú=13 \defÚ{\'U}
\catcode`\Ù=13 \defÙ{\`U}
\catcode`\Û=13 \defÛ{\^U}
\catcode`\Ü=13 \defÜ{\"U}

\catcode`\Ç=13\defÇ{\C C}

\begin{document}
\nocite{*}
\title[Transverse Poisson Structures]
  {Transverse Poisson structures to adjoint orbits  in semi-simple Lie algebras}
\author{Pantelis A.~Damianou}
\address{Pantelis Damianou, Department of Mathematics and Statistics\\
          University of Cyprus\\
          P.O.~Box 20537, 1678 Nicosia\\Cyprus}
\email{damianou@ucy.ac.cy}
\author{Hervé Sabourin}\author{Pol Vanhaecke}
\address{Hervé  Sabourin, Pol Vanhaecke,  Laboratoire de Mathématiques\\ UMR 6086 du CNRS\\
          Université de Poitiers
          \\86962 Futuroscope Chasseneuil Cedex\\France}
\email{sabourin@math.univ-poitiers.fr, pol.vanhaecke@math.univ-poitiers.fr}
\thanks{The authors would like to thank the Cyprus Research Foundation and the Minist\`ere Français des Affaires
\'etrang\`eres for their support.}
%
\subjclass[2000]{53D17, 17B10, 14J17}
\begin{abstract}
We study the transverse Poisson structure to adjoint orbits in a complex semi-simple Lie algebra. The problem is
first reduced to the case of nilpotent orbits. We prove then that in suitably chosen quasi-homogeneous coordinates
the quasi-degree of the transverse Poisson structure is $-2$. In the particular case of {\emph subregular}
nilpotent orbits we show that the structure may be computed by means of a simple determinantal formula, involving the
restriction of the Chevalley invariants on the slice. In addition, using results of Brieskorn and Slodowy, the
Poisson structure is reduced to a three dimensional Poisson bracket, intimately related to the simple rational
singularity that corresponds to the subregular orbit.
\end{abstract}
\date{\today}
\maketitle
\section{Introduction}\label{sec:intro}
The transverse Poisson structure was introduced by A.~Weinstein in \cite{WE}, stating in his famous splitting
theorem that every (real smooth or complex holomorphic) Poisson manifold $M$ is, in the neighbourhood of each point
$m$, the product of a symplectic manifold and a Poisson manifold of rank~$0$ at~$m$. The two factors of this
product can be geometrically realized as follows: let $S$ be the symplectic leaf through $m$ and let $N$ be any
submanifold of $M$ containing $m$ such that
$$
  T_m(M) = T_m(S) \oplus T_m(N).
$$
There exists a neighbourhood $V$ of $m$ in $N$, endowed with a Poisson structure, and a neighbourhood $U$ of $m$ in
$S$ such that, near $m$, $M$ is isomorphic to the product Poisson manifold $U \times V.$ The submanifold $N$ is
called \emph{a transverse slice} at $m$ to the symplectic leaf $S$ and the Poisson structure on $V\subset N$ is
called \emph{the transverse Poisson structure} to~$S$. Up to Poisson isomorphism, it is independent of the point
$m\in S$ and of the chosen transverse slice $N$ at $m$: given two points $m,\,m'\in S$ and two transverse slices
$N,\,N'$ at $m$ resp.\ $m'$ to $S$, there exist neighbourhoods $V$ of $m$ in $N$ and $V'$ of $m'$ in $N'$ such that
$(V,m)$ and $(V',m')$ are Poisson diffeomorphic.

\smallskip

When $M$ is the dual $\g^*$ of a complex Lie algebra $\g$, equipped with its standard Lie-Poisson structure, we
know that the symplectic leaf through $\mu\in\g^*$ is the co-adjoint orbit $\G\cdot\mu$ of the adjoint Lie group
$\G$ of $\g$. In this case, a natural transverse slice to $\G\cdot\mu$ is obtained in the following way: we choose
any complement $\n$ to the centralizer $\g(\mu)$ of $\mu$ in $\g$ and we take $N$ to be the affine subspace
$\mu+\n^\perp$ of $\g^*$. Since $\g(\mu)^\perp=\ad_\g^*\mu$ we have
$$
  T_\mu(\g^*) = T_\mu(\G\cdot\mu) \oplus T_\mu(N),
$$
so that $N$ is indeed a transverse slice to $\G\cdot \mu$ at $\mu$. Furthermore, defining on $\n^\perp$ any system
of linear coordinates $(q_1,\dots,q_k)$, and using the explicit formula for Dirac reduction (see Formula
(\ref{dirac_formula}) below), one can write down explicit formulas for the Poisson matrix $\Lambda_N:=
\left(\pb{q_i,q_j}_N\right)_{1\leq i,j\leq k}$ of the transverse Poisson structure, from which it follows easily that
the coefficients of $\Lambda_N$ are actually rational functions in $(q_1,\dots,q_k)$. As a corollary, in the
Lie-Poisson case, the transverse Poisson structure is always rational (\cite{SG}). One immediately wonders in which
cases the Poisson structure on $N$ is polynomial; more precisely, for which Lie algebras $\g$, for which co-adjoint
orbits for which complement $\n$.

\smallskip

Partial answers have been given in the literature for (co-) adjoint orbits in a semi-simple Lie algebra. In
(\cite{DA}), P.~Damianou computed explicitly the transverse Poisson structure to nilpotent orbits of $\gl_n$, for
$n\leq7$, corresponding to a particular complement $\n$, yielding that in this case the transverse Poisson
structure is polynomial. In (\cite{CU-RO}), R.~Cushman and M.~Roberts proved that there exists for any nilpotent
adjoint orbit of a semi-simple Lie algebra a special choice of a complement $\n$ such that the corresponding
transverse Poisson structure is polynomial. For the latter case, H.~Sabourin gave in (\cite{SA}) a more general
class of complements for which the transverse structure is polynomial, using in an essential way the machinery of
semi-simple Lie algebras; he also showed that the choice of complement $\n$ is relevant for the polynomial
character of the transverse Poisson structure by giving an example in which where the latter structure is rational
for a generic choice of complement.

\smallskip

When the transverse Poisson structure is polynomial one is tempted to define its degree as the maximal degree of
the coefficients $\pb{q_i,q_j}_N$ of its Poisson matrix, as was done in \cite{DA} and \cite{CU-RO}, who also
formulate several conjectures about this degree. Unfortunately, as shown in \cite{SA}, this degree depends strongly
on the choice of the complement $\n$, hence it is not intrinsically attached to the transverse Poisson
structure. We show in Section~\ref{sec:quasi} that the right approach is by using the more general notion of
``quasi-degree", i.e., we assign natural quasi-degrees $\varpi(q_i)$ to the variables $q_i$ ($i=1,\dots,k$) and we
show that, in the above mentioned class of complements, the quasi-degree of the transverse Poisson structure is
always $-2$, irrespective of the simple Lie algebra, the chosen adjoint orbit and the chosen transverse slice $N$!
In fact, the weights $\varpi(q_i)$ have a Lie-theoretic origin and are also independent of the particular
complement. It follows that $\pb{q_i,q_j}_N$ is a quasi-homogeneous polynomial of quasi-degree
$\varpi(q_i)+\varpi(q_j)-2$, for $1\leq i,j\leq k$.

\smallskip

Another result, established in this article, is that the study of the transverse Poisson structure to any adjoint
orbit $\G\cdot x$ can be reduced, via the Jordan decomposition of $x\in\g$, to the case of an adjoint
\emph{nilpotent} orbit. Thereby we explain why we are merely interested in the case of nilpotent orbits.




\smallskip

It is easy to observe that the transverse structure to the regular
nilpotent orbit $\O_{reg}$ of $\g$ is always trivial. So, the next
step is to consider the case of the \emph{subregular} nilpotent
orbit $\O_{sr}$ of $\g$. Then $N\cong\C^{\ell+2}$, where $\ell$ is
the rank of $\g$. The dimension of $\O_{sub}$ is two less than the
dimension of the regular orbit, so that the transverse Poisson
structure has rank $2$. It has $\ell$ independent polynomial
Casimirs functions $\chi_1, \dots, \chi_{\ell}$, where $\chi_i$ is
the restriction of the $i$-th Chevalley invariant $G_i$ to the
slice $N$.  In this case the transverse Poisson structure may be
obtained by a simple determinantal formula instead of the usual
rather complicated Dirac's constraint formula.  The determinantal
formula may be formulated as follows: in terms of linear
coordinates $q_1,\,q_2, \dots, q_{\ell+2}$ on $N$, the formula
\begin{equation}\label{eq:det_formula_intro}
 \pb{f,g}_{det}:= \frac{df \wedge dg \wedge d\chi_1 \wedge \dots \wedge d\chi_\ell}{dq_1\we dq_2\we\dots\we d
   q_{\ell+2}}
\end{equation}
defines a Poisson bracket on $N$, which coincides (up to a non-zero constant), with the transverse Poisson
structure on $N$.

\smallskip

As an application of Formula (\ref{eq:det_formula_intro}), we show in Theorem \ref{thm:3x3} that the Poisson matrix
of the transverse Poisson on $N$ takes, in suitable coordinates, the block form
  \begin{equation*}
    \widetilde\Lambda_N=
      \begin{pmatrix} 0 &0 \\ 0 &\Omega
      \end{pmatrix},\quad\hbox{where}\quad
      \renewcommand{\arraystretch}{2}
      \Omega=
      \begin{pmatrix}
        0 &\ds\pp{F}{q_{\ell+2}}&-\ds\pp{F}{q_{\ell+1}}\\
        -\ds\pp{F}{q_{\ell+2}}&0&\ds\pp{F}{q_\ell}\\
        \ds\pp{F}{q_{\ell+1}}&-\ds\pp{F}{q_\ell}&0
      \end{pmatrix}.
  \end{equation*}
The polynomial $F=F(u_1,\dots,u_{\ell-1},q_{\ell},q_{\ell+1},q_{\ell+2})$ that appears in this formula is precisely
the polynomial that describes the universal deformation of the (homogeneous or inhomogeneous) simple singularity of
the singular surface $N\cap\mathcal N$, where $u_1,\dots,u_{\ell-1}$ are the deformation parameters, which are
Casimirs for the Poisson structure on $N$, and where $\mathcal N$ is the nilpotent cone of $\g$. In particular, the
restriction of this Poisson structure to $N\cap\mathcal N$, is given by
\begin{equation*}
  \pb{x,y}=\pp {F_0}z,\qquad
  \pb{y,z}=\pp {F_0}x,\qquad
  \pb{z,x}=\pp {F_0}y,
\end{equation*}%
where $F_0(x,y,z):=F(0,\dots,0,x,y,z)$ is the polynomial that defines $N\cap\mathcal N$, as a surface in $\C^3$.
As we will recall in Paragraph \ref{sec:3x3}, Brieskorn \cite{brieskorn} showed that, in the ADE case, the
so-called adjoint quotient, $G=(G_1,\dots,G_\ell):\g\to\C^\ell$, restricted to the slice $N$, is a semi-universal
deformation of the singular surface $N\cap\mathcal N$; this result was generalized by Slodowy \cite{SL2} to the
other simple Lie algebras. Our Theorem \ref{thm:3x3} adds a Poisson dimension to this result.

\smallskip

The article is organized as follows: In Section \ref{sec:semi-simple}, we recall a few basic facts concerning
transverse Poisson structures and we show that the case of a general orbit in a semi-simple Lie algebra can be
reduced to the case of a nilpotent orbit.  In Section \ref{sec:quasi}, we recall the notion of quasi-homogeneity
and we show that, for a natural class of slices, the transverse Poisson structure is quasi-homogeneous of
quasi-degree $-2$.  In Section \ref{sec:examples} and the end of Section \ref{sec:det}, we show in some examples,
namely the Lie algebras $\g_2$, $\so_8$ and $\Liesl_4$, how the transverse Poisson structure can be computed
explicitly, and we use these examples to illustrate our results. In Section \ref{sec:det} we prove that in the case
of the subregular orbit, the transverse Poisson structure is given by a determinantal formula; we also show in that
section that this Poisson structure is entirely determined by the singular variety of nilpotent elements of the
slice.

\section{Transverse Poisson structures in semi-simple Lie algebras}\label{sec:semi-simple}
In this section we recall the main setup for studying the transverse Poisson structure to a (co-) adjoint orbit in
the case of a complex semi-simple Lie algebra $\g$, and we show how the case of a general orbit is related to the
case of a nilpotent orbit. We use the Killing form $\INN$ of $\g$ to identify $\g$ with its dual $\g^*$, which
leads to a Poisson structure on $\g$, given for functions $F,G$ on $\g$ at $x\in\g$ by
\begin{equation}\label{LP_gen_form}
  \pb{F,G}(x):=\inn x{\lb{dF(x),dG(x)}},
\end{equation}
where we think of $dF(x)$ and $dG(x)$ as elements of $\g\cong \g^*\cong T_x^*\g$. Since the Killing form is
$\Ad$-invariant, the isomorphism $\g\cong\g^*$ identifies the adjoint orbits $\G\cdot x$ of $\G$ with the
co-adjoint orbits $\G\cdot\mu$, so the symplectic leaf of $\PB$ that passes through~$x$ is the adjoint orbit
$\G\cdot x$. Also, as a transverse slice at $x$ to $\G\cdot x$ we can take an affine subspace $N:=x+\n^\perp$,
where $\n$ is any complementary subspace to the centralizer $\g(x):=\set{y\in\g\mid \lb{x,y}=0}$ of $x$ in $\g$ and
$\perp$ is the orthogonal complement with respect to the Killing form.  In order to give an explicit formula for
the transverse Poisson structure $\PB_{N}$ to $\G\cdot x$, let $(Z_1,\dots,Z_k)$ be a basis for $\g(x)$ and let
$(X_1,\dots,X_{2r})$ be a basis for $\n$, where $2r=\dim(\G\cdot x)$ is the rank of the Poisson structure
(\ref{LP_gen_form}) at $x$. These bases lead to linear coordinates $q_1,\dots,q_{k+2r}$ on $\g$, centered at $x$,
defined by $q_i(y):=\inn{y-x}{Z_i}$, for $i=1,\dots,k$ and $q_{k+i}(y):=\inn{y-x}{X_i}$, for $i=1,\dots,2r$. Since
$dq_i(y)=Z_i$ for $i=1,\dots,k$ and $dq_{k+i}(y)=X_i$ for $i=1,\dots,2r$, it follows from (\ref{LP_gen_form}) that
the Poisson matrix of $\PB$ at $y\in\g$ is given by
\begin{equation}\label{dirac_mat}
  \left(\pb{q_i,q_j}(y)\right)_{1\leq i,j\leq k+2r}=
  \begin{matr}{cc}
    A(y)&B(y)\\
    -\transp{B(y)}&C(y)
  \end{matr},
\end{equation}%
where
\begin{equation*}
  \renewcommand{\arraystretch}{1.4}
  \begin{array}{rcll}
     A_{i,j}(y)&=&\inn{y}{[Z_i,Z_j]}\qquad &1\leq i,j\leq k;\\
     B_{i,m}(y)&=&\inn{y}{[Z_i,X_m]} &1\leq i\leq k,\quad 1\leq m\leq 2r;\\
     C_{l,m}(y)&=&\inn{y}{[X_l,X_m]} &1\leq l,m\leq 2r.
  \end{array}
\end{equation*}%
It is easy to see that the skew-symmetric matrix $C(x)$ is invertible, so $C(y)$ is invertible for $y$ in a
neighborhood of $x$ in $\g$, and hence for $y$ in a neighborhood $V$ of $x$ in $N$.  By Dirac reduction, the
Poisson matrix of $\PB_N$ at $n\in V$, in terms of the coordinates $q_1,\dots,q_k$ (restricted to $V$), is given by
\begin{equation}\label{dirac_formula}
  \Lambda_N(n) = A(n) + B(n)C(n)^{-1}\transp{B(n)}.
\end{equation}
According to the Jordan-Chevalley decomposition theorem we can write $x$ as $x=s+e,$ where $s$ is semi-simple, $e$
is nilpotent and $[s,e] = 0$. Moreover, the respective centralizers of $x,s$ and $e$ are related as follows:
\begin{equation}\label{stab_rel}
  \g(x) = \g(s) \cap \g(e).
\end{equation}
This leads to a natural class of complements $\n$ to $\g(x)$; Since the restriction of $\INN$ to $\g(s)$ is
non-degenerate (see \cite[Prop. 1.7.7.]{DI}), we have a vector space decomposition of~$\g$ as
$$
  \g = \g(s) \oplus \n_s,
$$
where $\n_s = \g(s)^\perp$. Notice that $\n_s$ is $\g(s)$-invariant, $[\g(s),\n_s]\subset\n_s$, since
\begin{equation*}
  \inn{\g(s)}{[\g(s),\n_s]}=\inn{[\g(s),\g(s)]}{\n_s}\subset\inn{\g(s)}{\n_s}=\set{0}.
\end{equation*}%
Choosing any complement $\n_e$ of $\g(x)$ in $\g(s)$ we get the following decomposition of $\g$:
\begin{equation*}
  \g=\g(x)\oplus \n_e\oplus\n_s.
\end{equation*}%
We take then $\n:=\n_e\oplus\n_s$ and we denote $N_x:=x+\n^\perp$.  It follows that, if $n\in N_x$, so that
$n\in\g(s)$, then $\inn{n}{[\g(s),\n_s]}\subset\inn {\g(s)}{\n_s}=\set0,$ so that in particular
\begin{equation}\label{zero_12}
  \inn{n}{[\g(x),\n_s]}=\set0\quad\hbox{ and }\quad \inn{n}{[\n_e,\n_s]}=\set0.
\end{equation}%
Let us assume that the basis vectors $X_1,\dots,X_{2r}$ of $\n$ have been chosen such that
$X_1,\dots,X_{2p}\in\n_e$ and $X_{2p+1},\dots,X_{2r}\in\n_s$. Then the formulas (\ref{zero_12}) imply that the
Poisson matrix (\ref{dirac_mat}) takes at $n\in N_x$ the form
\begin{equation*}
  \Lambda(n)=
  \begin{matr}{ccc}
    A(n)&B_e(n)&0\\
    -\transp{B_e(n)}&C_{e}(n)&0\\
    0&0&C_s(n)
  \end{matr},
\end{equation*}%
where
\begin{equation*}
  \renewcommand{\arraystretch}{1.4}
  \begin{array}{rcll}
     A_{i,j}(n)&=&\inn{n}{[Z_i,Z_j]}\qquad &1\leq i,j\leq k;\\
     B_{e;i,m}(n)&=&\inn{n}{[Z_i,X_m]} &1\leq i\leq k,\quad 1\leq m\leq 2p;\\
     C_{e;l,m}(n)&=&\inn{n}{[X_l,X_m]} &1\leq l,m\leq 2p;\\
     C_{s;l,m}(n)&=&\inn{n}{[X_l,X_m]} &2p+1< l,m\leq 2r.
  \end{array}
\end{equation*}%
It follows from (\ref{dirac_formula}) that the Poisson matrix of the transverse Poisson structure on $N_x$ is given
by
\begin{equation}\label{dirac_formula_red}
  \Lambda_{N_x}(n) = A(n) + B_e(n)C_e(n)^{-1}\transp{B_e(n)}.
\end{equation}
Let us now restrict our attention to the Lie algebra $\g(s)$, which is reductive, so it decomposes as
$$
  \g(s) = {\mathfrak z}(s)\oplus\g_{ss}(s),
$$
where ${\mathfrak z}(s)$ is the center of $\g(s)$ and $\g_{ss}(s) = [\g(s),\g(s)]$ is the semi-simple part of
$\g(s)$. At the group level we have a similar decomposition of $\G(s)$, the centralizer of $s$ in $\G$, whose Lie
algebra is $\g(s)$, namely
$$
  \G(s) = \Z(s)\G_{ss}(s),
$$
%
where $\Z(s)$ is a central subgroup of $\G(s)$ and $\G_{ss}(s)$ is the semi-simple part of $\G(s)$, with Lie
algebra $\g_{ss}(s)$. Since $e \in \g(s)$ we can consider $\G(s)\cdot e$ as an adjoint orbit of the reductive Lie
algebra $\g(s)$. We may think of it as an adjoint orbit of a \emph{semi-simple} Lie algebra, since $\G(s)\cdot e =
\G_{ss}(s)\cdot e$; similarly we may think of a transverse slice to the adjoint orbit $\G(s)\cdot e$ as a
transverse slice to $\G_{ss}(s)\cdot e$, up to a summand with trivial Lie bracket. Denoting by $\perp_s$ the
orthogonal complement with respect to $\INN$, restricted to $\g(s)$, we have that $N:=e+\n_e^{\perp_s}$ is a
transverse slice to $\G(s)\cdot e$ since
$$
  \g(s)=\g(x)\oplus\n_e={\mathfrak z}(s) \oplus \g_{ss}(s)(e)\oplus \n_e. 
$$
We have used that $\g(x)=\g(s)(e)$, the centralizer of $e$ in $\g(s)$, which follows from (\ref{stab_rel}). In
terms of the bases $(Z_1,\dots,Z_k)$ of $\g(x)$ and $(X_1,\dots,X_{2p})$ of $\n_e$ that we had picked, the Poisson
matrix takes at $n\in N$ the form
\begin{equation*}
  \begin{matr}{cc}
    A(n)&B_e(n)\\
    -\transp{B_e(n)}&C_{e}(n)\\
  \end{matr},
\end{equation*}%
which leads by Dirac reduction to the following formula for the transverse Poisson structure $\Lambda_{N}$ on $N$,
$$
  \Lambda_N(n) = A(n) + B_e(n)C_e(n)^{-1}\transp{B_e(n)},
$$
where $n\in N$. This yields formally the same formula as (\ref{dirac_formula_red}), except that it is evaluated
at points $n$ of $N$, rather than at points of $N_x$. However, since $\n_e^{\perp_s} = \g(s) \cap \n_e^\perp =
\n_s^\perp \cap \n_e^\perp = (\n_s + \n_e)^\perp = \n^\perp$, the affine subspaces $N_x$ and $N$ only differ by a
translation, $N_x=s+e+\n^\perp=s+N$, so that they, and their Poisson matrices with respect to the coordinates
$q_1,\dots,q_k$, can be identified. It leads to the following proposition.
\begin{proposition}\label{prp:gen_to_nilp}
  Let $x\in\g$ be any element, $\G\cdot x$ its adjoint orbit and $x=s+e$ its Jordan-Chevalley decomposition. Given
  any complement $\n_e$ of $\g(x)$ in $\g(s)$ and putting $\n:= \n_s \oplus \n_e$, where $\n_s=\g(s)^\perp$, the
  parallel affine spaces $N_x:= x + \n^\perp$ and $N:= e+\n^\perp$ are respectively transverse slices to the
  adjoint orbit $\G\cdot x$ in $\g$ and to the nilpotent orbit $\G(s)\cdot e$ in $\g(s)$. The Poisson structure on
  both transverse slices is given by the same Poisson matrix, namely that of (\ref{dirac_formula_red}), in terms of
  the same affine coordinates restricted to the corresponding transverse slice.
\end{proposition}

\smallskip In short, the transverse Poisson structure to any adjoint orbit $\G\cdot x$ of a semi-simple (or
reductive) algebra $\g$ is essentially determined by the transverse Poisson structure to the underlying nilpotent
orbit $\G(s)\cdot e$ defined by the Jordan decomposition $x = s+e$. A refinement of this proposition will be given in
Corollary \ref{cor:gen_to_nilp}.

\section{The polynomial and the quasi-homogeneous character of the tranverse Poisson structure}\label{sec:quasi}
In this section we show that, for a natural class of transverse slices to a nilpotent orbit $\O$, which we
equip with an adapted set of linear coordinates, centered at a nilpotent element $e\in\O$, the transverse Poisson
structure is quasi-homogeneous (of quasi-degree~$-2$), in the following sense.
\begin{definition}
  Let $\nu =(\nu_1,\dots, \nu_d)$ be non-negative integers. A polynomial $P\in\C[x_1,\dots,x_d]$
  is said to be \emph{quasi-homogeneous} (relative to $\nu$) if for some integer $\kappa$,
  $$
    \forall t \in {\C}, P(t^{\nu_1}x_1, \dots, t^{\nu_d}x_d) = t^\kappa P(x_1,\dots,x_d).
  $$
  The integer $\kappa$ is then called the \emph{quasi-degree} (relative to $\nu$) of $P$, denoted
  $\varpi(P)$. Similarly, a polynomial Poisson structure $\PB$ on $\C[x_1,\dots,x_d]$ is said to be
  \emph{quasi-homogeneous} (relative to $\nu$) if there exists $\kappa\in\Z$ such that, for any quasi-homogeneous
  polynomials $F$ and $G$, their Poisson bracket $\pb{F,G}$ is quasi-homogeneous of degree
\begin{equation*}
  \wght(\pb{F,G})=\wght(F)+\wght(G)+\kappa;
\end{equation*}%
  equivalently, for any $i,j$ the polynomial $\pb{x_i,x_j}$ is quasi-homogeneous of quasi-degree
  $\nu_i+\nu_j+\kappa$. Then $\kappa$ is called the \defi{quasi-degree} of $\PB$.
\end{definition}
We first show that, given $\O$, we can choose a system of linear coordinates on $\g$, centered at some nilpotent
element $e\in\O$, such that the Lie-Poisson structure on $\g$ is quasi-homogeneous relative to some vector $\nu$
that has a natural Lie-theoretic interpretation. In order to describe how this happens we need to recall a few
facts on the theory of semi-simple Lie algebras that will be used throughout this paper.  First, one chooses a
Cartan subalgebra $\h$ of $\g$, with corresponding root system $\Delta({\h})$, from which a basis $\Pi({\h})$ of
simple roots is selected.  The \emph{rank} of $\g$, which is the dimension of $\h$, is denoted by $\ell$.
According to the Jacobson-Morosov-Kostant correspondence (see \cite[pars.\ 32.1 and 32.4]{TA-YU}), there is a
canonical triple $(h,e,f)$ of elements of $\g$, associated with $\O$ and completely determined by the following
properties:
\begin{enumerate}
  \item[(1)] $(h,e,f)$ is a $\Liesl_2$-triple, i.e., $[h,e]=2e$, $[h,f]=-2f$ and $[e,f]=h$;
  \item[(2)] $h$ is the characteristic of $\O$, i.e., $h\in\h$ and $\a(h) \in \{0,1,2\}$ for any simple root $\a\in\Pi(\h)$.
  \item[(3)] $\O = \G\cdot e$.
\end{enumerate}
The triple $(h,e,f)$ leads to two decompositions of $\g$:

\smallskip
(1) A decomposition of $\g$ into eigenspaces relative to $\ad_h$. Each eigenvalue being an integer we have
  $$
    \g = \bigoplus_{i \in {\Z}} \g(i),
  $$
  where $\g(i)$ is the eigenspace of $\ad_h$ that corresponds to the eigenvalue $i$. For example, $e\in\g(2)$ and
  $f\in\g(-2)$.
\smallskip

(2) Let $\s$ be the Lie subalgebra of $\g$ isomorphic to $\Liesl_2$, which is generated by $h,e$ and~$f$.  The Lie
  algebra $\g$ is an $\s$-module, hence it decomposes  as
  $$
    \g = \bigoplus_{j=1}^k V_{n_j},
  $$
  where each $V_{n_j}$ is a simple $\s$-module, with $n_j+1=\dim V_{n_j}$ $\ad_h$-weights $n_j, n_j-2,
  n_j-4,\dots,-n_j$. Moreover, $k=\dim\g(e)$, since the centralizer $\g(e)$ is generated by the highest weight
  vectors of each $V_{n_j}$. It follows that
  \begin{equation}\label{weight_to_dim_orbit}
     \sum_{j=1}^{k} n_j = \dim \g - k=\dim (\G\cdot e)=2r.
  \end{equation}
We pick a system of linear coordinates on $\g$, centered at $e$, by using Slodowy's action. We recall the
construction of this action from \cite{SL}.  First, he considers the one-parameter subgroup of $\G$,
\begin{equation*}
  \begin{array}{lcccl}
    \l&:&\C^*&\to&\G\\
    & &t &\mapsto&\exp(\l_th)
  \end{array}
\end{equation*}%
where $\lambda_t$ is a complex number such that $e^{-\lambda_t} = t$. The restriction of $\Ad$ to this subgroup
leaves every eigenspace $\g(i)$ invariant, and acts for each $t$ as a homothecy with ratio $t^{-i}$ on $\g(i)$:
\begin{equation}\label{homothecy}
  \forall x \in \g(i),\ \Ad_{\lambda(t)} x=t^{-i}x.
\end{equation}%
Since $e\in\g(2)$ the action $\rho$ of $\C^*$ on $\g$, defined for $t\in\C^*$ and for $y\in\g$ by $\rho_t\cdot
y:=t^{2}\Ad_{\lambda(t)}y$ fixes~$e$; we refer to $\rho$ as \emph{Slodowy's action}. In order to see how it leads
to quasi-homogeneous coordinates, let us denote for $x\in\g$ by $\FF_x$, the function defined by
$\FF_x(y):=\inn{y-e}x$, for $y\in\g$. Then (\ref{homothecy}) and $\Ad$-invariance of the Killing form imply that if
$x\in\g(i)$ then
\begin{eqnarray*}
  \left(\rho_t^*\FF_x\right)(y)&=&\inn{\rho_{t^{-1}}\cdot y-e}x=t^{-2}\inn{\Ad_{\l(t^{-1})}(y-e)}x\\
   &=&t^{-2}\inn{y-e}{\Ad_{\l(t)}x}=t^{-2}\inn{y-e}{t^{-i}x}=t^{-i-2}\FF_x(y).
\end{eqnarray*}
It follows that the quasi-degree $\wght(\FF_x)$ of $\FF_x$ is $i+2$, for $x\in\g(i)$. According to
(\ref{LP_gen_form}), one has, for any $x,\,y,\,z\in\g$,
\begin{equation}\label{eq:pb_for_Fx}
  \pb{\FF_x,\FF_y}(z)=\inn{z}{\lb{x,y}}=\FF_{[x,y]}(z)+\inn{e}{[x,y]}.
\end{equation}%
If $x\in\g(i)$ and $y\in\g(j)$, with $i+j\neq -2$ then $\inn{e}{[x,y]}=0$ and so
\begin{eqnarray*}
  \wght(\pb{\FF_x,\FF_y})-\wght(\FF_x)-\wght(\FF_y)&=&\wght(\FF_{\lb{x,y}})-\wght(\FF_x)-\wght(\FF_y)\\
    &=&i+j+2-(i+2)-(j+2)=-2.
\end{eqnarray*}%
This result extends to the case $i+j=-2$, since then $\wght(\FF_{[x,y]})=i+j+2=0$, which is the quasi-degree of the
constant function $\inn{e}{[x,y]}$. This proves the following proposition.
\begin{proposition}\label{prp:quasi_lie}
  Let $\g$ be a semi-simple Lie algebra, identified with its dual using its Killing form, let $\O$ be a nilpotent
  adjoint orbit of $\g$, with canonical triple  $(h,e,f)$. Let $x_1,\dots,x_d$ be any basis of $\g$, where each $x_k$
  belongs to some eigenspace $\g(i_k)$ of $\ad_h$ and let $\FF_k$ be the dual coordinates on $\g$, centered at $e$,
  $\FF_k(y):=\inn{y-e}{x_k}.$ Then the Lie-Poisson structure $\PB$ on $\g$ is quasi-homogeneous of degree $-2$ with
  respect to $(\wght(\FF_1),\dots,\wght(\FF_d))=(i_1+2,\dots,i_d+2)$.
\end{proposition}
\qed

We now wish to show that, upon picking a suitable transverse slice $N$ to $\O$ at $e$, the transverse Poisson
structure on $N$ is also quasi-homogeneous (of degree $-2$). Following \cite{SA} we consider the set ${\mathcal
N}_h$ of all subspaces $\n$ of $\g$ that are complementary to $\g(e)$ in $\g$, and which are $\ad_h$-invariant. For
$\n\in{\mathcal N}_h$ we let $N:=e+\n^\perp$, which is a transverse slice to $\G\cdot e$. The $\ad_h$-invariance of
$\n$ implies on the one hand that $\rho$ leaves $N$ invariant: if $y\in e+\n^\perp$ then
\begin{equation*}
  0=\inn{y-e}{\Ad_{\l(t^{-1})}\n}=\inn{\Ad_{\l(t)}(y-e)}{\n}=t^{-2}\inn{\rho_t\cdot y-e}{\n},
\end{equation*}%
so that indeed $\rho_t\cdot y\in e+\n^\perp$. On the other hand, it implies that $\n$ admits a basis where each
basis vector belongs to an eigenspace of $\h$. Thus we can specialize the above basis $x_1,\dots,x_d$ so that it be
adapted to $\n$: we can choose a basis $(Z_1,\dots,Z_k)$ for $\g(e)$ and a basis $(X_1, \dots, X_{2r})$ for $\n$ in
such a way that:
\begin{enumerate}
  \item[(1)] each $Z_i, 1 \leq i \leq k,$ is a highest weight vector of weight $n_i$ ;
  \item[(2)] each $X_i, 1 \leq i \leq 2r$, is a weight vector of weight $\nu_i$.
\end{enumerate}
The linear coordinates (centered at $e$) $\FF_{Z_1},\dots,\FF_{Z_k}$, restricted to $N$, will be denoted by
$q_1,\dots,q_k$. In view of the above, their quasi-degrees are defined as $\wght(q_i):=n_i+2$. The fact that the
transverse Poisson structure is polynomial in terms of these coordinates was first shown in \cite[Thm 2.3]{SA}. In
the following proposition we give a refinement of this statement.

\begin{proposition}\label{prp:atp_poly}
  In the notation of Proposition \ref{prp:quasi_lie}, the transverse Poisson structure on $N:=e+\n^\perp$, where
  $\n\in{\mathcal N}$, is a polynomial Poisson structure that is quasi-homogeneous of degree $-2$, with respect to
  the quasi-degrees $n_1+2,\dots,n_k+2$, where $n_1,\dots,n_k$ denote the highest weights of $\g$ as an
  $\s$-module.
\end{proposition}
\begin{proof}
According to (\ref{dirac_formula}) we need to show that for any $1\leq i,j\leq k$ the functions $A_{ij}$ and
$(BC^{-1}\transp{B})_{ij}$ are quasi-homogeneous of degree $\wght(q_i)+\wght(q_j)-2=n_i+n_j+2$. For $A_{ij}$ this
is clear, since $A$ is part of the Poisson matrix of the Lie-Poisson structure on $\g$, which we have seen to be
quasi-homogeneous of degree $-2$. Similarly, we have that $\wght(B_{ip})= n_i+\nu_p+2$. Since
$$
  \wght(B_{ip}C^{-1}_{ps}B_{js})=n_i + n_j+\nu_p+\nu_s+4+\wght(C^{-1}_{ps}),
$$
it means that we need to show that
\begin{equation}\label{C_wght}
  \wght(C^{-1}_{ps})=-\nu_p-\nu_s-2.
\end{equation}%
This follows from the fact that $\sum_{i=1}^{2r}(\nu_i+1)=0$, itself a consequence of
(\ref{weight_to_dim_orbit}). Indeed, consider a term of the form $C'_{ij} = C_{i_1j_1}\dots C_{i_{2r-1}j_{2r-1}}$,
where
\begin{equation*}
  \set{i_1,i_2,\dots,i_{2r-1}}=\set{1,2,\dots,2r}\setminus\set s\hbox{ and }
  \set{j_1,i_2,\dots,j_{2r-1}}=\set{1,2,\dots,2r}\setminus\set p.
\end{equation*}%
Then
\begin{equation*}
  \wght(C'_{ij})=\sum_{k=1}^{2r-1}(\nu_{i_k} + \nu_{j_k} +2) = -\nu_s -\nu_p -2,
\end{equation*}%
A typical term of $C^{-1}_{ps}$ is of the form $\frac{C'_{ij}}{\Delta(C)}$, where $\Delta(C)$ is the determinant of
$C$ which is constant, as it is of quasi-degree zero, by the same argument (this observation was already made in
\cite[Thm 2.3]{SA}). This gives us (\ref{C_wght}).
\end{proof}

\smallskip

Let us consider now the case of any adjoint orbit $\G\cdot x$ and $x = s+e$ the Jordan-Chevalley decomposition of
$x$, a case that we already considered in Proposition \ref{prp:gen_to_nilp}. A well-known result
(\cite[par. 32.1.7.]{TA-YU}) says that there exists a $\Liesl_2$-triple $(h,e,f)$ such that $[s,h] = [s,f] =
0$. Consequently, $(h,e,f)$ is a $\Liesl_2$-triple of the reductive Lie algebra $\g(s)$ and we can also suppose
that, up to conjugation by elements of $\G(s)$, $h$ is the characteristic of $\G(s)\cdot e$. Let ${\mathcal
N}_{s,h}$ the set of all complementary subspaces to $\g(x)$ in $\g(s)$ which are $\ad_h$-invariant.  Then, by
applying Proposition \ref{prp:atp_poly}, we get the following result.

\begin{corollary}\label{cor:gen_to_nilp}
  As in Proposition \ref{prp:gen_to_nilp}, let $\n_s = \g(s)^\perp, \n_e \in {\mathcal N}_{s,h}$ and $ \n = \n_s
  \oplus \n_e$. Let $N_x:= x + \n^\perp$, which is a transverse slice to $\G\cdot x$. Then, the transverse Poisson
  structure on $N_x$ is polynomial and is quasi-homogeneous of quasi-degree $-2$.
\end{corollary}

From now, a transverse Poisson structure given by Proposition \ref{prp:atp_poly} will be called an {\it adjoint
transverse Poisson structure} or ATP-structure.

\section{Examples}\label{sec:examples}
In this section we want to show in two examples how one computes the ATP-structure in practice. In the first
example we consider the subregular orbit of $\g_2$, and we do the computations without chosing a representation of
$\g_2$. In the second example, the subregular orbit of $\so_8$, we use a concrete representation, rather than
referring to tables for the explicit formulas of the Lie brackets in a Chevalley basis. These two examples will
also serve later in the paper as an illustration of the results that we will prove on the nature of the
ATP-structure. Both examples correspond to subregular orbits, which lead to two of the simplest non-trivial
ATP-structures, in the following sense. If $\O$ is an adjoint orbit in $\g$ then the ATP-structure to $\O$ has rank
$\dim\g-\ell-\dim\O$ at a generic point of any transverse slice to $\O$, since the Lie-Poisson structure on $\g$
has rank\footnote{Recall that $\ell$ denotes the rank of $\g$.}  $\dim\g-\ell$, at a generic point of $\g$. For the
regular nilpotent orbit $\O_{reg}$, the ATP-structure is trivial, because $\dim\O_{reg}=\dim\g-\ell$. So, the first
interesting nilpotent orbit to consider is the subregular orbit, denoted by $\O_{sr}$. We recall two well-known
facts (see~\cite{CO-MC}):
%
\begin{enumerate}
  \item[-] the subregular orbit $\O_{sr}$ is the unique nilpotent orbit which is open and dense in the complement of
  $\O_{reg}$ in the nilpotent cone;
  \item[-] $\dim \O_{sr} = \dim\g-\ell-2$.
\end{enumerate}
It follows that the ATP-structure of the subregular orbit is of dimension $\ell+2$ and its generic rank is $2$.  In
each of the two examples that follow we give the characteristic triplet $(h,e,f)$ that corresponds to the orbit, we
derive from it a basis of the $\ad_h$-weight spaces, which leads to basis vectors $Z_i$ of $\g(e)$ and $X_j$ of an
$\ad_h$-invariant complement to $\g(e)$ in $\g$. The Lie brackets of these elements then lead to the matrices
$A,\,B$ and $C$ in (\ref{dirac_mat}), which by Dirac's formula (\ref{dirac_formula}) yields the matrix $\Lambda_N$
of the transverse Poisson structure.

\subsection{The subregular orbit of type $G_2$}
We first consider the case of the subregular orbit of the Lie algebra $\g:=\g_2$. Denoting the basis of simple roots by
$\Pi=\set{\a,\b}$, where $\b$ is the longer root, its Dynkin diagram is given by
%

\smallskip
\begin{center}
  \epsfbox{g2.fig}
\end{center}

\noindent
and it has the following positive roots : $\Delta_+ = \{\a,\,\b,\,\a+\b,\,2\a+\b,\,3\a+\b,\,3\a+2\b,\,3\a+\b,$
$3\a+2\b\}$. The vectors that constitute the Chevalley basis\footnote{The choice of Chevalley basis that we use is
explicitly described in \cite[Chapter VII.4]{TA}.} of $\g$ are denoted by $H_{\a},\,H_{\b}$ for the Cartan
subalgebra, $X_\gamma$ for the 6 positive roots $\gamma\in\Delta_+$ and $Y_\gamma$ for the six negative roots
$-\gamma\in\Delta_+$. According to \cite[Chap. 8.4]{CO-MC}, the characteristic $h$ of the subregular orbit
$\O_{sr}$ is given by the sequence of weights $(0,2)$, which means that $\can{\a,h}=0$ and $\can{\b,h}=2$, which
yields $h=2H_\a+4H_\b$. The decomposition of $\g$ into $\ad_h$-weight spaces $\g(i)$ consists of the following five
subspaces
\begin{align}\label{eq:G2_bases}
\g(4)& = \langle X_{3\a+ 2\b}\rangle,\nonumber\\
\g(2)& = \langle X_{\b},\, X_{\a+\b},\,X_{2\a+\b},\, X_{3\a+\b}\rangle,\nonumber\\
\g(0)& = \langle H_\a,\,H_\b,\,X_{\a},\, Y_{\a}\rangle,\\
\g(-2)& = \langle Y_{\b},\, Y_{\a+\b},\,Y_{2\a+\b},\, Y_{3\a+\b}\rangle,\nonumber\\
\g(-4)& = \langle Y_{3\a+ 2\b}\rangle.\nonumber
\end{align}
Taking for $e$ and $f$ an arbitrary linear combination of the above basis elements of $\g(2)$, resp.\ of $\g(-2)$
and expressing that $[e,f]=h$ one easily finds that the $\Liesl_2$-triple corresponding to $\O_{sr}$ is given by
$$
  e = X_{\b} + X_{3\a + \b},\ h = 2H_{\a} + 4H_{\b},\ f = 2Y_{\b}+2Y_{3\a+\b}.
$$
Picking the vectors in the positive subspaces $\g(i)$ that commute with $e$ leads to the following basis vectors of
$\g(e)$:
\begin{align}\label{eq:G2_bases_Z}
  Z_1 &= X_{\b} +  X_{3\a + \b}.\nonumber\\
  Z_2 &= X_{2\a +\b}, \nonumber\\
  Z_3 &= X_{\a+ \b}, \\
  Z_4 &= X_{3\a + 2\b}, \nonumber
\end{align}
We obtain an $\ad_h$-invariant complementary subspace $\n$ of $\g(e)$ by completing these vectors with additional
vectors that are taken from the bases (\ref{eq:G2_bases}) of the subspaces $\g(i)$. Our choice of basis vectors for
$\n$, ordered by weight, is as follows:
\begin{align*}
  X_1 &= X_{\b}, &X_{6} &= Y_{\b},\\
  X_2 &= X_{\a}, &  X_{7} &= Y_{\a+\b},\\
  X_3 &= H_{\a}, &X_{8} &= Y_{2\a+\b}, \\
  X_4 &= H_{\b}, &X_{9} &= Y_{3\a+\b}, \\
  X_{5} &= Y_{\a}, &X_{10} &= Y_{3\a+2\b}.
\end{align*}
The Lie brackets of these basis vectors for $\g$, which are listed in \cite[Chap. VII.4]{TA}, yield the Poisson
matrix
$ \begin{matr}{cc} A&B\\ -\transp{B}&C\end{matr}$
of the Lie-Poisson structure on $\g$, in terms of the coordinates
$\FF_{Z_1},\dots,\FF_{Z_4},\FF_{X_1},\dots,\FF_{X_{10}}$ on $\g$, as
$A_{ij}=\pb{\FF_{Z_i},\FF_{Z_j}}=\FF_{\lb{Z_i,Z_j}}+\inn{e}{[Z_i,Z_j]}$ (see (\ref{eq:pb_for_Fx})), and similarly
for the other elements of the Poisson matrix. We give the restriction of the matrices $A,\,B$ and $C$ to the
transverse slice $N:=e+\n^\perp$ only, which amounts to keeping in the Lie brackets only the vectors
$Z_1,\dots,Z_4$, as $\FF_{X}(n)=\inn{e-n}{X}=0$ for $X\in\n$ and $n\in N=e+\n^\perp$. In terms of the coordinates
$q_1,\dots,q_4$ on $N$, where $q_i$ is the restriction of $\FF_{Z_i}$ to $N$, we get
\begin{equation*}
  A=\begin{matr}{cccc}
    0&0&0&0\\
    0&0&-3q_4&0\\
    0&3q_4&0&0\\
    0&0&0&0
  \end{matr},
\end{equation*}%
\begin{equation*}
  B=\begin{matr}{cccccccccc}
    0&0&0&-q_4&0&q_1&-q_2&q_3&0&0\\
    0&3q_1&-q_2&0&2q_3&0&0&0&0&0\\
    0&2q_2&q_3&-q_3&0&0&0&0&0&0\\
    q_4&q_3&-3q_1&q_1&q_2&0&0&0&0&0
  \end{matr},
\end{equation*}%
\begin{equation*}
  C=\frac13\begin{matr}{cccccccccc}
    0&3q_3&0&0&0&0&0&0&0&1\\
    -3q_3&0&0&0&0&0&3&0&0&0\\
    0&0&0&0&0&3&0&0&-3&0\\
    0&0&0&0&0&-2&0&0&1&0\\
    0&0&0&0&0&0&0&3&0&0\\
    0&0&-3&2&0&0&0&0&0&0\\
    0&-3&0&0&0&0&0&0&0&0\\
    0&0&0&0&-3&0&0&0&0&0\\
    0&0&3&-1&0&0&0&0&0&0\\
    -1&0&0&0&0&0&0&0&0&0
  \end{matr}.
\end{equation*}%
%
%
Substituted in (\ref{dirac_formula}) this yields the following  Poisson matrix for the ATP-structure:
\begin{equation}\label{eq:g2_atp}
  \Lambda_N =
  \begin{pmatrix}
    0 & 0 & 0&0 \\
    0 &0 &-3q_4 &2 q_1 q_2-2 q_3^2 \\
    0 &3 q_4&0 & 2 q_2^2 -2 q_1 q_3\\
    0& -2 q_1 q_2+2 q_3^2& -2 q_2^2+2 q_1 q_3 &0
  \end{pmatrix}.
\end{equation}
It follows from (\ref{eq:G2_bases}) and (\ref{eq:G2_bases_Z}) that the quasi-degree of $q_1,\, q_2$ and $q_3$ is
$4$, while the quasi-degree of $q_4$ is $6$. One easily reads off from (\ref{eq:g2_atp}) that, with respect to
these quasi-degrees, the ATP-structure is quasi-homogeneous of quasi-degree $-2$.
\subsection{The subregular orbit of type $D_4$}
We take now $\g=\so_8$ and we realize $\g$ as the following set of matrices :
\begin{equation*}
  \left\{\begin{pmatrix}
    Z_1 &Z_2 \\ Z_3 & -\transp{Z_1}
  \end{pmatrix}
  \mid Z_i \in \text{Mat}_4({\C}),\text{ with } Z_2,Z_3 \text{ skew-symmetric} \right\}.
\end{equation*}
Let $\h$ denote the Cartan subalgebra of $\g$ that consists of all diagonal matrices in $\g$. It is clear that $\h$
is spanned by the four matrices $H_i:= E_{i,i} - E_{4+i,4+i}, 1 \leq i \leq 4.$ Define for $i=1,\dots 4$ the linear
map $e_i\in\h^*$ by $e_i(\sum a_kH_k) = a_i$. Then the root system of $\g$ is:
$$
  \Delta:= \{ \pm e_i \pm e_j \mid 1 \leq i,j \leq 4, i \not= j\},
$$
and a basis of simple roots is $\Pi:=\{\a_1,\a_2,\a_3,\a_4\}$, where
$$
  \a_1 = e_1 - e_2,\ \a_2 = e_2 - e_3,\ \a_3 = e_3 - e_4,\ \a_4 = e_3 + e_4.
$$
It leads to the following  Chevalley basis of $\g$:
\begin{align*}
  X_{e_i -e_j} &= E_{i,j} - E_{4+j,4+i}, \\
  X_{e_i +e_j} &= E_{i,4+j} - E_{j,4+i}, \qquad i<j\\
  X_{-e_i -e_j} &= -E_{4+i,j} + E_{4+j,i}, \qquad i<j \\
  H_{e_i -e_j} &= H_i - H_j, \\
  H_{e_i +e_j} &= H_i + H_j.
\end{align*}

According to (\cite[Chap. 5.4]{CO-MC}), the characteristic $h$ of the subregular orbit is given by the sequence
of weights $(2,0,2,2)$. It follows that
$$
  h = 4H_{\alpha_1}+ 6H_{\alpha_2} + 4H_{\alpha_3} + 4H_{\alpha_4}.
$$
The positive $\ad_h$-weight spaces are
\begin{align}\label{eq:d4_bases}
  \g(0)& = {\h} \oplus \langle X_{\a_2}, X_{-\a_2}\rangle,\nonumber\\
  \g(2)& = \langle X_{\a_1}, X_{\a_3},X_{\a_4},X_{\a_1+ \a_2}, X_{\a_2 + \a_3},X_{\a_2+\a_4}\rangle,\nonumber\\
  \g(4)& = \langle X_{\a_1+ \a_2+\a_3}, X_{\a_2 +\a_3+\a_4},X_{\a_1+ \a_2+\a_4}\rangle,\\
  \g(6)& = \langle X_{\a_1+ \a_2+\a_3+\a_4}, X_{\a_1 +2 \a_2+\a_3+\a_4} \rangle.\nonumber
\end{align}
As in the first example it follows that the following is the  canonical $\Liesl_2$-triple, associated to $\O_{sr}$:
\begin{eqnarray*}
  e &=& X_{\a_1} + X_{\a_1 + \a_2} - X_{\a_2 + \a_4}+ 2X_{\a_3} - X_{\a_4},\\
  h &=& 4H_{\a_1}+ 6H_{\a_2} + 4H_{\a_3} + 4H_{\a_4},\\
  f &=& X_{-\a_1} + 3X_{-\a_1 - \a_2} - 3X_{-\a_2 -\a_4}+ 2X_{-\a_3} - X_{-\a_4}.
\end{eqnarray*}
We can now define the basis vectors $Z_i$ of $\g(e)$ and $X_j$ of an $\ad_h$-invariant complementary subspace $\n$
to $\g(e)$, in terms of the Chevalley basis, as follows:

\begin{align}\label{eq:d4_bases_Z}
  Z_1 &= X_{\a_1 +\a_2} - X_{\a_2+\a_4}+ 2X_{\a_3}, \nonumber\\
  Z_2 &=  X_{\a_1 + \a_2 + \a_3 + \a_4}, \nonumber\\
  Z_3 &= X_{\a_1 + 2\a_2 + \a_3 + \a_4}, \nonumber\\
  Z_4 &= X_{\a_1} - X_{\a_4}, \\
  Z_5 &= X_{\a_2+ \a_3}+  X_{\a_2+ \a_4}-  X_{\a_3}- X_{\a_4},\nonumber \\
  Z_6 &= X_{\a_1 + \a_2 + \a_3 } +  X_{\a_1 + \a_2 +\a_4} -  X_{\a_2 + \a_3 + \a_4}, \nonumber
\end{align}
\begin{align*}
  X_1 &= X_{\a_1 + \a_2 + \a_3}, &X_{12} &= X_{-\a_1},\\
  X_2 &= X_{\a_2 + \a_3 + \a_4}, &  X_{13} &= X_{-\a_3},\\
  X_3 &= X_{\a_4} &  X_{14}, &= -X_{-\a_4},\\
  X_4 &= X_{\a_3}, &X_{15} &= X_{-\a_1-\a_2},\\
  X_5 &= X_{\a_2 + \a_4}, & X_{16} &= X_{-\a_2-\a_3},\\
  X_6 &= H_{\a_1}, &X_{17} &= -X_{-\a_2-\a_4},\\
  X_7 &= H_{\a_2}, & X_{18} &= -X_{-\a_1-\a_2-\a_3},\\
  X_8 &= H_{\a_3}, & X_{19} &= -X_{-\a_1 -\a_2-\a_4},\\
  X_9 &= H_{\a_4}, & X_{20} &= -X_{-\a_2 -\a_3-\a_4},\\
  X_{10} &= X_{\a_2}, &X_{21} &= -X_{-\a_1 -\a_2-\a_3-\a_4},\\
  X_{11} &= X_{-\a_2}, &X_{22} &= -X_{-\a_1 -2\a_2-\a_3-\a_4}.
\end{align*}
If we denote by $\bar Z_1,\dots,\bar Z_6$ the dual basis (w.r.t.\ $\inn XY=\frac12\Trace(XY)$) of the basis
$Z_1,\dots,Z_6$ of $\g(e)$ then a typical element of the transverse slice $N = e +\n^\perp$ is given by
$e+\sum_{i=1}^6q_i\bar Z_i$, i.e.,
\begin{equation}\label{matrix_Q}
  Q =
  \begin{pmatrix}
    0 &1 &1 &0 &0 &0 &0 &0 \\
    q_4 &0 &0 &0 &0 &0 &0 &-1 \\
    q_1 &0 &0 &2 &0 &0 &0 &-1 \\
    0 &q_5 &0 &0 &0 &1 &1 &0 \\
    0 &-q_3 &-q_2 &0 &0 &-q_4 &-q_1 &0 \\
    q_3 &0 &q_6 &0 &-1 &0 &0 &-q_5 \\
    q_2 &-q_6 &0 &0 &-1 &0 &0 &0 \\
    0 &0 &0 &0 &0 &0 &-2 &0
  \end{pmatrix}
\end{equation}
and we can compute the matrix $A$, restricted to $N$, by $A_{ij}=\inn{Q}{[Z_i,Z_j]}$, and similarly for the
matrices $B$ and $C$. A direct substitution in (\ref{dirac_formula}) leads to the following Poisson matrix for the
the ATP-structure:
\begin{equation}\label{eq:atp_matrix_d4}
  \Lambda_N =\frac12
  \begin{pmatrix}
    0 &  q_4 q_6 &- q_4 q_6& 0& -2 q_6 &2q_{16} \\
    - q_4 q_6  &0 &0 &  q_4 q_6& -q_5 q_6 &-2 q_{36} \\
     q_4 q_6 &0 &0 & - q_4 q_6&  q_5 q_6  & 2  q_{36}\\
    0& - q_4 q_6&  q_4 q_6&0 &2  q_6&  -2 q_{16}\\
     2 q_6 & q_5 q_6&- q_5q_6& -2 q_6&0&2 q_{56}\\
    -2 q_{16} &2q_{36}&-2 q_{36}&2q_{16}&-2q_{56} &0
  \end{pmatrix},
\end{equation}
where
\begin{eqnarray}\label{q1j_def}
 q_{16} &=& 2q_2 - q_1q_4 - q_4q_5 + q_4^2,\nonumber\\
 q_{36} &=& q_3q_4 -q_2q_4-q_2q_5,\\
 q_{56} &=& 2q_3 -2q_2 -q_5^2+q_4q_5 - q_1q_5.\nonumber
\end{eqnarray}
It follows from (\ref{eq:d4_bases}) and (\ref{eq:d4_bases_Z}) that the quasi-degrees of the variables $q_i$ are:
$\wght(q_1)=\wght(q_4)=\wght(q_5)=4,$ $\wght(q_2)=\wght(q_3)=8$ and $\wght(q_6)=6$. The fact that the ATP-structure
is quasi-homogeneous of quasi-degree $-2$ is again easily read off from (\ref{eq:atp_matrix_d4}).
\section{The subregular case}\label{sec:det}
In this section we will give an explicit description of the ATP-structure in the case of the subregular orbit
$\O_{sr}\subset\g$, where $\g$ is a semi-simple Lie algebra. Since in the case of the subregular orbit the generic
rank of the ATP-structure on the transverse slice $N$ is two and since we know $\dim(N)-2$ independent Casimirs,
namely the basic $\Ad$-invariant functions on $\g$, restricted to $N$, we will easily derive that the ATP-structure
is the determinantal structure (also called Nambu structure), determined by these Casimirs, up to multiplication by
a \emph{function.}  What is much less trivial to show is that this function is actually just a \emph{constant.} For
this we will use Brieskorn's theory of simple singularities, which is recalled in Paragraph \ref{par:singularities}
below. First we recall the basic facts on $\Ad$-invariant functions on $\g$ and link them to the ATP-structure.

\subsection{Invariant functions and Casimirs}
Let $\O_{sr} = \G\cdot e$, be a subregular orbit in the semi-simple Lie algebra $\g$, let $(h,e,f)$ be the
corresponding canonical $\Liesl_2$-triple and consider the transverse slice $N := e + {\n}^\perp$ to $\G\cdot e$,
where $\n$ is an $\ad_h$-invariant complement to $\g(e)$. We know from Section \ref{sec:quasi} that the
ATP-structure on $N$, equipped with the linear coordinates $q_1,\dots,q_{k}$, is a quasi-homogeneous polynomial
Poisson structure of generic rank $2$.  Let $S(\g^*)^\G$ be the algebra of $\Ad$-invariant polynomial functions on
$\g$. By a classical theorem due to Chevalley, $S(\g^*)^\G$ is a polynomial algebra, generated by $\ell$
homogeneous polynomials $(G_1,\dots,G_\ell)$, whose degree $d_i:=\deg(G_i)=m_i+1$, where $m_1,\dots,m_\ell$ are the
exponents of $\g$. These functions are Casimirs of the Lie-Poisson structure on $\g$, since $\Ad$-invariance of
$G_i$ implies that $[x,dG_i(x)]=0$ and hence the Lie-Poisson bracket (\ref{LP_gen_form}) is given by
\begin{equation*}
  \pb{F,G_i}(x)=\inn x{\lb{dF(x),dG_i(x)}}=-\inn{\lb{x,dG_i(x)}}{dF(x)}=0,
\end{equation*}%
for any function $F$ on $\g$. If we denote by $\chi_i$ the restriction of $G_i$ to the transverse slice $N$ then it
follows that these functions are Casimirs of the ATP-structure. The polynomials $\chi_i$ are not homogeneous, but
they are quasi-homogeneous, as shown in the following proposition.
\begin{lemma}\label{lma:qd_chi}
  Each $\chi_i$ is a quasi-homogeneous polynomial of quasi-degree $2d_i$, relative to the quasi-degrees $(2+n_1,\dots,
  2+n_k)$.
\end{lemma}
\begin{proof}
Since $\chi_i$ is of degree $d_i$ and $\chi_i$ is $\Ad$-invariant, we get
$$ \rho_t^*(\chi_i) = \chi_i\circ\rho_{t^{-1}}=\chi_i\circ(t^{-2}\Ad_{\l^{-1}(t)})=t^{-2d_i}
  \chi_i\circ\Ad_{\l^{-1}(t)}=t^{-2d_i}\chi_i,
$$
so that $\chi_i$ has quasi-degree $2d_i$.
\end{proof}

\subsection{Simple singularities}\label{par:singularities}
Let $\h$ be a Cartan subalgebra of $\g$. The Weyl group $\W$ acts on $\h$ and the algebra $S(\g^*)^\G$ of
$\Ad$-invariant polynomial functions on $\g$ is isomorphic to $S(\h^*)^\W$, the algebra of $\W$-invariant
polynomial functions on $\h^*$. The inclusion homomorphism $S(\g^*)^\G\hookrightarrow S(\g^*)$, is dual to a morphism
$\g\to\h/\W$, called the \emph{adjoint quotient}. Concretely, the adjoint quotient is given by
\begin{equation}\label{eq:adjoint_quotient}
  \begin{array}{lcccl}
    G&:&\g&\to&\C^\ell\\
    & &x &\mapsto&(G_1(x),\,G_2(x),\dots,G_\ell(x)).
  \end{array}
\end{equation}%
%
The zero-fiber $G^\mi(0)$ of $G$ is exactly the nilpotent variety ${\mathcal N}$ of $\g$. As we are interested in
$N\cap {\mathcal N}=N\cap G^\mi(0)=\chi^{-1}(0)$, which is an affine surface with an isolated, simple singularity,
let us recall the notion of simple singularity (see \cite{SL2} for details). Up to conjugacy, there are five types
of finite subgroups of $\SL2=\SL2({\C})$, denoted by ${\mathcal C}_p, {\mathcal D}_p, {\mathcal T}, {\mathcal O}$
and ${\mathcal I}$. Given such a subgroup $\F$, one looks at the corresponding ring of invariant polynomials
${\C}[U,V]^\F$. In each of the five cases, $\C[u,v]^\F$ is generated by three fundamental polynomials $X,Y,Z$,
subject to only one relation $R(X,Y,Z) = 0$, hence the quotient space $\C^2/\F$ can be identified, as an affine
surface, with the singular surface in $\C^3$, defined by $R=0$. The origin is its only singular point; it is called
a \emph{(homogeneous) simple singularity}. The exceptional divisor of the minimal resolution of $\C^2/\F$ is a
finite set of projective lines. If two of these lines meet, then they meet in a single point, and
transversally. Moreover, the intersection pattern of these lines forms a graph that coincides with one of the
simply laced Dynkin diagrams of type $A_\ell,\,D_\ell,\,E_6,\,E_7$ or $E_8$. This type is called the type of the
singularity. Moreover, every such Dynkin diagram (i.e., of type ADE) appears in this way; see Table
\ref{sing_table_1}.
\def\mytabvrule{\vrule height0.6cm depth0.45cm width0pt}
\begin{table}[ht]
  \caption{The basic correspondence between finite subgroups $\F$ of $\SL2$, homogeneous simple singularities,
  defined by an equation $R(X,Y,Z)=0$ and simply laced simple Lie algebras of type $\Delta$.\label{sing_table_1}}
  \begin{center}
  \begin{tabular}{|c|c|c|}
    \hline
      \ Group $\F$\ \ &\ Singularity $R(X,Y,Z) = 0$\ \ &\ Type $\Delta$\ \ \mytabvrule \\
    \hline
      ${\mathcal C}_{\ell+1}$ &$X^{\ell+1} +YZ = 0$ &$A_{\ell}$\mytabvrule\\
      ${\mathcal D}_{\ell-2}$ &$X^{\ell-1} +XY^2 + Z^2 = 0$ &$D_{\ell}$\mytabvrule\\
      ${\mathcal T}$ &$X^4 +Y^3 + Z^2 = 0$ &$E_6$\mytabvrule\\
      ${\mathcal O}$ &$X^3Y +Y^3 + Z^2 = 0$ &$E_7$\mytabvrule\\
      ${\mathcal I}$ &$X^5 +Y^3 + Z^2 = 0$  & $E_8$\mytabvrule\\
    \hline
  \end{tabular}
  \end{center}
\end{table}

\smallskip

For the other simple Lie algebras (of type $B_\ell,\, C_\ell,\, F_4$ or $G_2$), there exists a similar
correspondence. By definition, an \emph{(inhomogeneous) simple singularity} of type $\Delta$ is a couple
$(V,\Gamma)$ consisting of a homogeneous simple singularity $V = {\C}^2/\F$ and a group $\Gamma = \F'/\F$ of
automorphisms of $V$ according to Table \ref{inhomogeneous_table}.

\begin{table}[ht]
  \caption{We give the list of all possible inhomogeneous singularities of type $\Delta=(V,\Gamma)$, where $V$ is
   one of the homogeneous simple singularities and $\Gamma=\F'/\F$ is a group of automorphisms of $V$. The labels
   $B_\ell,\,C_\ell,\,F_4$ and $G_2$ for these types will become clear in Proposition
   \ref{prp:simple_sing}. \label{inhomogeneous_table}}
  \begin{center}
  \begin{tabular}{|c|c|c|c|c|}
    \hline
     Type $\Delta$ \ & $V$ &$\F$&$\F'$&$\Gamma= \F'/\F$\mytabvrule \\
    \hline
      $B_\ell$ &$A_{2\ell-1}$&${\mathcal C}_{2\ell}$&${\mathcal D}_\ell$&${\Z}/2{\Z}$\mytabvrule \\
      \ $C_\ell$\ \ &\ $D_{\ell+1}$\ \  &\ ${\mathcal D}_{\ell-1}$\ \ &\ ${\mathcal D}_{2\ell-2}$\ \ &\
                ${\Z}/2{\Z}$\ \ \mytabvrule \\
      $F_4$&$E_6$& ${\mathcal T}$&${\mathcal O}$&${\Z}/2{\Z}$\mytabvrule\\
      $G_2$&$D_4$&${\mathcal D}_2$&${\mathcal O}$&$\Z/3\Z$\mytabvrule\\
    \hline
  \end{tabular}
  \end{center}
\end{table}
The connection between the diagram of $(V,\Gamma)$ and that of $V$ can be described as follows: the action of
$\Gamma$ on $V$ lifts to an action on a minimal resolution of $V$ which permutes the components of the exceptional
set. Then, we obtain the diagram of $(V,\Gamma)$ as a $\Gamma$-quotient of that of $V$. It leads to Table
\ref{sing_table_2}, which is the non-simply laced analog of Table \ref{sing_table_1}.

\begin{table}[ht]
  \caption{For each of the inhomogeneous simple singularities, of type $\Delta$ (see Table
  \ref{inhomogeneous_table}), the correspondending homogeneous simple singularity $V=\C^2/\F$ is given by its
  equation $R(X,Y,Z)=0$, together with the action of $\Gamma=\F'/\F$ on $V$. In the last line, $\a$ is a
  non-trivial cubic root of unity.\label{sing_table_2}}
  \begin{center}
  \begin{tabular}{|c|c|c|}
    \hline
      Type $\Delta$ \ & Singularity $R(X,Y,Z) = 0$&$\Gamma$-action\mytabvrule\\
    \hline
       $B_\ell$&$X^{2\ell} +YZ = 0$&$(X,Y,Z) \lra (-X,Z,Y)$\mytabvrule\\
       $C_\ell$&$X^\ell +XY^2 +Z^2 = 0$&$(X,Y,Z) \lra (X,-Y,-Z)$\mytabvrule\\
       $F_4$&$X^4 + Y^3 + Z^2 = 0$&$(X,Y,Z) \lra (-X,Y,-Z)$\mytabvrule\\
       $G_2$&$X^3 + Y^3 + Z^2 = 0$&$(X,Y,Z) \lra (\a X,\a^2 Y,Z)$\mytabvrule \\
    \hline
  \end{tabular}
  \end{center}
\end{table}

We can now state the following extension of a theorem of Brieskorn, which is due to Slodowy (\cite[Thms 1 and
2]{SL2}).

\begin{proposition}\label{prp:simple_sing}
  Let $\g$ be a simple complex Lie algebra, with Dynkin diagram of type~$\Delta$. Let $\O_{sr}= \G\cdot e$ be the
  subregular orbit and $N = e + {\mathfrak n}^\perp$ a transverse slice to $\G\cdot e$.  The surface $N \cap
  {\mathcal N} = \chi^{-1}(0)$ has a (homogeneous or inhomogeneous) simple singularity of type $\Delta$.
\end{proposition}
To finish this paragraph, we illustrate the above results in the case of the two examples that were given in
Section \ref{sec:examples}. In both cases we give the invariants, restricted to the slice~$N$, and their zero
locus, the surface $\chi^{-1}(0)$.

\smallskip
1) For the subregular orbit of $\g_2$, the invariant functions, restricted to the slice~$N$, are given by
\begin{eqnarray}\label{eq:g2_cas}
  \chi_1 &=& q_1, \nonumber\\
  \chi_2 &=& 12q_1q_2q_3-4q_2^3-4q_3^3+9q_4^2,
\end{eqnarray}
which leads to an affine surface $\chi^{-1}(0)$ in $\C^4$, isomorphic to the surface in $\C^3$ defined by
$$
  4q_2^3+4q_3^3-9q_4^2=0.
$$
Up to a rescaling of the coordinates, this is the polynomial $R$ that was given in Table \ref{sing_table_2}.

%

\smallskip

2) For the subregular orbit of $\so_8$, the invariant functions, restricted to the slice~$N$, are found as the
   (non-constant) coefficients of the characteristic polynomial of the matrix $Q$ (see (\ref{matrix_Q})):
\begin{eqnarray}\label{eq:d4_cas}
  \chi_1 &=& -2q_1 -2q_4,\nonumber\\
  \chi_2 &=& -12q_2 -4q_3 -4q_4q_5 + (q_1+q_4)^2,\nonumber\\
  \chi_3 &=& -q_2+q_3 -q_4q_5, \\
  \chi_4 &=& -4q_1q_2 -16q_2q_5 -12q_3q_4 + 12 q_2q_4 +4q_1q_3+ 4q_4^2q_5 + 4q_1q_4q_5 -4q_6^2.\nonumber
\end{eqnarray}
By linearly eliminating the variables $q_1,\,q_2$ and $q_3$ from the equations $\chi_i=0$, for $i=1,2,3$, we find
that $\chi^{-1}(0)$ is isomorphic to the affine surface in $\C^3$, defined by
$$
   4q_4^2q_5 - 2q_4 q_5^2 + q_6^2= 0.
$$
Its defining polynomial corresponds to the polynomial $R$ in Table \ref{sing_table_1}, upon putting $X=i\gamma
q_4,\, Y=\gamma(q_5-q_4)$ and $Z=q_6$, where $\gamma$ is any cubic root of $2i$.
%



\subsection{The determinantal Poisson structure}
We prove in this paragraph the announced result that the ATP-structure in the subregular case is a determinantal
Poisson structure, determined by the Casimirs. Let us first point out how such a determinantal Poisson structure is
defined. Let $C_1,\dots,C_{d-2}$ be $d-2$ (algebraically) independent polynomials in $d>2$ variables
$x_1,\dots,x_{d}$. For such a polynomial $F$, let us denote by $\nabla F$ its differential $\diff F$, expressed
in the natural basis $\diff x_i$, i.e., $\nabla F$ is a column vector with elements $(\nabla F)_i=\frac{\p F}{\p
x_i}$. Then a polynomial Poisson structure is defined on $\C^{d}$ by
\begin{equation}\label{det_poi_def}
  \pb{F,G}_{det}:=\det(\nabla F,\ \nabla G,\ \nabla C_1,\ \dots,\ \nabla C_{d-2}),
\end{equation}%
where $F$ and $G$ are arbitrary polynomials. It is clear that each of the $C_i$ is a Casimir of $\PB_{det}$, so
that in particular the generic rank of $\PB_{det}$ is two. Notice also that if the Casimirs $C_i$ are
quasi-homogeneous, with respect to the weights $\varpi_i:=\varpi(x_i)$ then for any quasi-homogeneous elements $F$
and $G$ we have that
\begin{equation*}
  \wght(\pb{F,G}_{det})=\wght(F)+\wght(G)+\sum_{i=1}^{d-2}\wght (C_i)-\sum_{i=1}^{d}\wght_i,
\end{equation*}%
an easy consequence of the definition of a determinant and of the following obvious fact: if $F$ is any
quasi-homogeneous polynomial, then $\frac{\p F}{\p x_i}$ is quasi-homogeneous and $\wght(\frac{\p F}{\p x_i}) =
\wght(F) - \wght_i$.  It follows that $\PB_{det}$ is quasi-homogeneous of quasi-degree $\kappa$, where
\begin{equation}\label{for:det_wght}
  \kappa=\sum_{i=1}^{d-2}\wght (C_i)-\sum_{i=1}^{d}\wght_i.
\end{equation}%
\smallskip

Applied to our case it means that we have two polynomial Poisson structures on the transverse slice $N$ which have
$G_1,\dots,G_{\ell}$ as Casimirs on $N\cong\C^{\ell+2},$ namely the ATP-structure and the determinantal structure,
constructed by using these Casimirs.

\smallskip

In our two examples (see Section \ref{sec:examples}), these structures are easily compared by explicit
computation. For the subregular orbit of $\g_2$, we have according to (\ref{det_poi_def}) that
\begin{equation*}
  \left(\Lambda_{det}\right)_{ij}=\det\left(\nabla q_i\ \nabla q_j\ \nabla\chi_1\ \nabla\chi_2\right),
\end{equation*}%
where $\chi_1$ and $\chi_2$ are the Casimirs (\ref{eq:g2_cas}). This leads to
$$
  \Lambda_{det} = -6
  \begin{pmatrix}
    0 & 0 & 0&0 \\
    0 &0 &-3q_4 &2 q_1 q_2-2 q_3^2 \\
    0 &3 q_4&0 & 2 q_2^2 -2 q_1 q_3\\
    0& -2 q_1 q_2+2 q_3^2& -2 q_2^2+2 q_1 q_3 &0
  \end{pmatrix}.
$$
In view of (\ref{eq:g2_atp}), it follows that $\Lambda_{det} = -6\Lambda_N$, so that both Poisson structures
coincide. For $\so_8$ one finds similarly, using the Casimirs $\chi_1,\dots,\chi_4$ in (\ref{eq:d4_cas})
$$
 \Lambda_{det} = -128
 \begin{pmatrix}
    0 &  q_4 q_6 &- q_4 q_6& 0& -2 q_6 &2q_{16} \\
    - q_4 q_6  &0 &0 &  q_4 q_6& -q_5 q_6 &-2 q_{36} \\
     q_4 q_6 &0 &0 & - q_4 q_6&  q_5 q_6  & 2  q_{36}\\
    0& - q_4 q_6&  q_4 q_6&0 &2  q_6&  -2 q_{16}\\
     2 q_6 & q_5 q_6&- q_5q_6& -2 q_6&0&2 q_{56}\\
    -2 q_{16} &2q_{36}&-2 q_{36}&2q_{16}&-2q_{56} &0
 \end{pmatrix},
$$
where $q_{16},\,q_{36}$ and $q_{56}$ are given by (\ref{q1j_def}). In view of (\ref{eq:atp_matrix_d4}), both
Poisson structures again coincide, $\Lambda_{det} = -256\Lambda_N$.

\smallskip

In order to show that, in the subregular case, the ATP-structure and the determinantal structure always coincide,
i.e., differ by a constant factor only, we first show, that both structures coincide, up to factor, which is a
rational function.
\begin{proposition}\label{prp:ss_form}
  Let $\PB$ and $\PB'$ be two non-trivial polynomial Poisson structures on $\C^{d}$ which have $d-2$ common
  independent polynomial Casimirs $C_1,\dots,C_{d-2}$. There exists a rational function $R\in\C(x_1,\dots,x_{d})$
  such that $\PB=R\PB'$.
\end{proposition}
\begin{proof}
  Let $M$ and $M'$ denote the Poisson matrices that correspond to $\PB$ and $\PB'$ with respect to the coordinates
  $x_1,\dots,x_{d}$. If we denote $\RR:=\C(x_1,\dots,x_{d})$ then $M$ and $M'$ both act naturally as a
  skew-symmetric endomorphism on the $\RR$-vector space~$\RR^d$. The subspace $H$ of $\RR^d$ that is spanned by
  $\nabla C_1,\dots,\nabla C_{d-2}$ is the kernel of both maps, hence we have two induced skew-symmetric
  endomorphisms $\varphi$ and $\varphi'$ of the quotient space $\RR^d/H$. Since the latter is two-dimensional,
  $\varphi'$ and $\varphi$ are proportional, $\varphi'=R\varphi$, where $R\in\RR$. Since $M$ and $M'$ have the same
  kernel, $M'=RM$.
\end{proof}
Applied to our two Poisson structures $\PB_N$ and $\PB_{det}$ the proposition yields that $\PB_N=R\PB_{det}$, where
$R=P/Q\in\RR$. We show in the following result that $R$ is actually a (non-zero) constant, thereby characterizing
completely the ATP-structure in the subregular case.

\begin{thm}\label{thm:determinantal}
  Let $\O_{sr}$ be the subregular nilpotent adjoint orbit of a complex semi-simple Lie algebra $\g$ and let
  $(h,e,f)$ be the canonical triple, associated to $\O_{sr}$. Let $N = e + {\n}^\perp$ be a transverse slice to
  $\O_{sr}$, where $\n$ is an $\ad_h$-invariant complementary subspace to $\g(e)$. Let $\PB_N$ and $\PB_{det}$
  denote respectively the ATP-structure and the determinantal structure on $N$. Then $\PB_N= c\PB_{det}$ for some
  $c \in {\C^*}$.
\end{thm}

\begin{proof}
By the above, $\PB_N=R\PB_{det}$, where $R\in\RR$. If $R$ has a non-trivial denominator $Q$, then all elements of the
Poisson matrix of $\PB_{det}$ must be divisible by $Q$, since both Poisson structures are polynomial. Then along
the hypersurface $Q=0$ the rank of $(\nabla \chi_1,\dots,\nabla \chi_\ell)$ is smaller than $\ell$, hence
$\chi^{-1}(0)$ is singular along the curve $\chi^{-1}(0)\cap (Q=0)$. However, by Proposition \ref{prp:simple_sing},
we know that $\chi^{-1}(0)$ has an isolated  singularity, leading to a contradiction. This shows that $Q$ is
a constant, and hence that $R$ is a polynomial.

In order to show that the polynomial $R$ is constant it suffices to show that the quasi-degrees of $\PB_N$ and
$\PB_{det}$ are the same, which amounts (in view of Proposition \ref{prp:quasi_lie}) to showing that the
quasi-degree of $\PB_{det}$ is $-2$. This follows from the following formula, due to Kostant (see \cite[Thm
7]{KO}), which expresses the dimension of the regular orbit in terms of the exponents $m_i$ of $\g$:
\begin{equation}
  2\sum_{i=1}^\ell m_i = \dim \O_{reg} = \dim \g - \ell.
\end{equation}
Indeed, if we use this formula, Lemma \ref{lma:qd_chi} and (\ref{weight_to_dim_orbit}) in the formula
(\ref{for:det_wght}) for the quasi-degree of $\PB_{det}$, then we find
\begin{eqnarray*}
  \kappa&=&\sum_{i=1}^\ell\varpi(\chi_i)-\sum_{i=1}^{\ell+2}\varpi(q_i)=2\sum_{i=1}^\ell d_i-\sum_{i=1}^{\ell+2}(n_i+2)\\
         &=&2\sum_{i=1}^\ell m_i-\sum_{i=1}^{\ell+2}n_i-4\\
         &=&\dim\g-\ell-(\dim\g-\ell-2)-4=-2.
\end{eqnarray*}
\end{proof}

\subsection {Reduction to a $3 \times 3$ Poisson matrix}\label{sec:3x3}

Let $\O_{sr}$ be the subregular nilpotent adjoint orbit of a complex semi-simple Lie algebra $\g$ of rank
$\ell$. Let $(h,e,f)$ be its associated canonical $\Liesl_2$-triple, and let $N:= e + {\n}^\perp$ be a transverse
slice to $\O_{sr}$, where $\n$ is an $\ad_h$-invariant complementary subspace to $\g(e)$. Let $\PB_N$ be the
ATP-structure defined on $N$. Recall that $N$ is equipped with linear coordinates $q_1,\dots,q_{\ell+2}$ defined in
Section \ref{sec:semi-simple}, and that $\PB_N$ has independent Casimirs $\chi_1, \dots, \chi_{\ell}$, which are
the restrictions to $N$ of the basic homogeneous invariant polynomial functions on $\g$.

\smallskip

Our goal now is to show that, in well-chosen coordinates, the ATP-structure $\PB_N$ on $N$ is essentially given by
a $3 \times 3$ skew-symmetric matrix, closely related to the polynomial that defines the singularity. More
precisely, we have the following theorem.

\begin{thm}\label{thm:3x3}
  After possibly relabeling the coordinates $q_i$ and the Casimirs $\chi_i$, the $\ell+2$ functions
  \begin{equation}\label{eq:new_coo}
    \chi_i, 1 \leq i \leq \ell-1,\quad\hbox{and}\quad q_\ell,\, q_{\ell+1},\, q_{\ell+2}
  \end{equation}
  form a system of (global) coordinates on the affine space $N$. The Poisson matrix of the ATP-structure on $N$
  takes, in terms of these coordinates, the form
  \begin{equation}\label{eq:omega}
    \widetilde\Lambda_N=
      \begin{pmatrix} 0 &0 \\ 0 &\Omega
      \end{pmatrix},\quad\hbox{where}\quad
      \renewcommand{\arraystretch}{2}
      \Omega=c'
      \begin{pmatrix}
        0 &\ds\pp{\chi_\ell}{q_{\ell+2}}&-\ds\pp{\chi_\ell}{q_{\ell+1}}\\
        -\ds\pp{\chi_\ell}{q_{\ell+2}}&0&\ds\pp{\chi_\ell}{q_\ell}\\
        \ds\pp{\chi_\ell}{q_{\ell+1}}&-\ds\pp{\chi_\ell}{q_\ell}&0
      \end{pmatrix},
  \end{equation}
  for some non-zero constant $c'$. It has the polynomial $\chi_\ell$ as Casimir, which reduces to the polynomial which
  defines the singularity, when setting $\chi_j=0$ for $j=1,2, \dots, \ell-1$.
\end{thm}
\begin{proof}
The non-Poisson part of this theorem is due to Brieskorn and Slodowy. Before proving the Poisson part of the
theorem, namely that the Poisson matrix takes the above form (\ref{eq:omega}), we explain, for the convenience of
the reader, the basic notions of singularity theory that are used in their proof, see \cite{SL2} for details. Let
$(X_0, x) $ be the germ of an analytic variety $X_0$ at the point $x$.  A deformation of $(X_0,x)$ is a pair
$(\Phi, \imath)$ where $\Phi:X\to U$ is a flat morphism of varieties, with $\Phi(x)=u$, and where the map
$\imath:X_0\to\Phi^{-1}(u)$ is an isomorphism.  Such a deformation is called \defi{semi-universal} if any other
deformation of $(X_0,x)$ is isomorphic to a deformation, induced from $(\Phi,\imath)$ by a local change of
variables, in a neighborhood of $x$. The semi-universal deformation of $(X_0,x)$ is unique up to isomorphism. It
can be explicitly described in the following case: let $(X_0,0)$ be a germ of a hypersurface of $\C^d$, which is
singular at $0$, say $X_0$ is locally given by $f(z)=0$. Then the semi-universal deformation of $(X_0,0)$ is the
(germ at the origin of the) map
\begin{equation}
  \begin{array}{lcccl}
    \Phi&:&\C^k\times\C^d&\to&\C^k\times\C\\
    & &(u,z) &\mapsto&(u,F(u,z)),
  \end{array}
\end{equation}
where
\begin{equation}\label{eq:F}
  F(u,z)=f(z)+ \sum_{i=1}^k g_i(z) u_i,
\end{equation}
and where the polynomials $1, g_1, g_2, \dots, g_k$ represent a vector space basis of the Milnor (or Tjurina) algebra
\begin{equation}\label{eq:milnor_algebra}
  \M(f):=\frac{\C[z_1, \dots, z_d]}{\left( f, \pp f{z_1}, \dots, \pp f{z_d} \right)}
        =\frac{\C[z_1, \dots, z_d]}{\left(\pp f{z_1}, \dots, \pp f{z_d} \right)},
\end{equation}
where the last equality is valid whenever $f$ is quasi-homogeneous, which is true in the present context. The
dimension $\dim\M(f)=k+1$ is called the \defi{Milnor number} of $f$.

We can now formulate Brieskorn's result. It says that the map $\chi:N\to\C^\ell$, which is the restriction of the
adjoint quotient (\ref{eq:adjoint_quotient}) to the slice $N$, is a semi-universal deformation of the singular
surface $N\cap{\mathcal N}$. More precisely, when the Lie algebra is of the type ADE, then the map
\begin{equation}
  \begin{array}{lcccl}
    \Phi&:&\C^{\ell-1}\times\C^3&\to&\C^{\ell-1}\times\C\\
    & &\left((\chi_1,\dots,\chi_{\ell-1}),(q_{\ell},q_{\ell+1},q_{\ell+2})\right)
                &\mapsto&\left( (\chi_1,\dots,\chi_{\ell-1}),\chi_\ell\right)
  \end{array}
\end{equation}
is the semi-universal deformation of the singular surface $N\cap {\mathcal N}$; for the other types one has to
consider $\Gamma$-invariant semi-universal deformations, as was shown by Slodowy, see Table
\ref{inhomogeneous_table} and \cite{SL2}. It is implicit in Brieskorn's statement that
$(\chi_1,\dots,\chi_{\ell-1},q_{\ell},q_{\ell+1},q_{\ell+2})$ form a system of coordinates on $N$, which comes from
the fact that one can solve the $\ell-1$ equations $\chi_i=\chi_i(q)$ \emph{linearly} for $\ell-1$ of the variables
$q_i$, namely the Casimirs have the form
\begin{equation}\label{eq:coo_change}
  \begin{pmatrix}
     \chi_1\\ \vdots \\ \chi_{\ell-1}
  \end{pmatrix}
  =A
  \begin{pmatrix}
    q_1\\ \vdots \\ q_{\ell-1}
  \end{pmatrix}
  +
  \begin{pmatrix}
    F_1(q_\ell,q_{\ell+1},q_{\ell+2})\\ \vdots \\ F_{\ell-1}(q_\ell,q_{\ell+1},q_{\ell+2})
  \end{pmatrix},
\end{equation}%
where $A$ is a constant matrix, with $\det A\neq0$; this will be illustrated in the examples below.

We now get to the Poisson part of the proof. Since the coordinate functions $\chi_1,\dots,\chi_{\ell-1}$ are
Casimirs, the Poisson matrix $\widetilde\Lambda_N$ takes with respect to these coordinates the block form
\begin{equation*}
  \widetilde\Lambda_N= \begin{pmatrix} 0 &0 \\ 0 & \Omega \end{pmatrix}
\end{equation*}
where $\Omega$ is a $3\times 3$ skew-symmetric matrix. We know from Theorem \ref{thm:determinantal} that the
ATP-structure is a constant multiple of the determinantal structure. Since $\det A\in\C^*$ it follows from
(\ref{eq:coo_change}) that for $\ell\leq i,j\leq \ell+2$,
\begin{equation*}
  \widetilde\Lambda_{ij}:=c\det(\nabla q_i\ \nabla q_j\ \nabla\chi_1\ \dots\ \nabla\chi_\ell)
       =c'\det(\nabla' q_i\ \nabla' q_j\ \nabla'\chi_{\ell}),
\end{equation*}%
where $c$ and $c'$ are non-zero constants and $\nabla'$ denotes the restriction of $\nabla$ to $\C^3$, namely
\begin{equation*}
  \nabla'F=
  \begin{pmatrix}
    \ds\pp F{q_{\ell}}&
    \ds\pp F{q_{\ell+1}}&
    \ds\pp F{q_{\ell+2}}
  \end{pmatrix}^\top.
\end{equation*}
The explicit formula (\ref{eq:omega}) for $\Omega$ follows from it at once.
\end{proof}

\noindent{\bf Examples}

\smallskip
\noindent
1) For the subregular orbit of $\g_2$ we have, according to  (\ref{eq:g2_cas}), that $\chi_1=q_1$, so that
   $\chi_2$, expressed in terms of $q_2,\,q_3,\,q_4$ and $\chi_1$, is given by
\begin{equation*}
  \chi_2=9 q_4^2-4q_2^3-4q_3^3+12  \chi_1 q_2 q_3.
\end{equation*}
The Poisson matrix (\ref{eq:g2_atp}) of the ATP-structure is already in the form (\ref{eq:omega}), with $c'=-1/6$
(and $\chi_1=q_1$).  Since the Milnor algebra (\ref{eq:milnor_algebra}) is given in this case by
$\M(9q_4^2-4q_2^3-4q_3^3)= \C[q_2,q_3,q_4]/\left(q_2^2,q_3^2,q_4\right)$, one easily sees that $1$ and the
coefficient $q_2 q_3$ of $u_1$ form indeed a vector space basis for the $\Gamma$-invariant elements of the Milnor
algebra (see Table \ref{sing_table_2}); cfr.\ \cite[p.\ 136]{SL2}.

\smallskip

\noindent 2) We now turn to the case of $\so_8$. Recall from (\ref{eq:d4_cas}) that its ATP structure has Casimirs
$\chi_1, \dots, \chi_4$. As stated in the proof of Theorem \ref{thm:3x3}, we can solve three of them linearly for
$q_1$, $q_2$, $q_3$ in terms of $\chi_1$, $\chi_2$, $\chi_3$ and the last three variables $q_4$, $q_5$, and
$q_6$. We obtain
\begin{eqnarray*}
  q_1&=& -q_4-\frac{\chi_1}{2}, \\
  q_2&=& \frac{1}{64} \left( \chi_1^2-16 \chi_3-4 \chi_2- 32 q_4 q_5 \right),  \\
  q_3&=& \frac{1} {64} \left( \chi_1^2+48 \chi_3-4 \chi_2 +32 q_4 q_5\right). \\
\end{eqnarray*}
Substituted in  $\chi_4$ this yields
\begin{equation*}
  \chi_4=8 q_4 q_5^2-16 q_4^2 q_5-4q_6^2-4 \chi_1q_4 q_5 +(\chi_2-\frac{1}{4} \chi_1^2 + 4\chi_3 )q_5 -16 \chi_3 q_4
   -2 \chi_1 \chi_3,
\end{equation*}
so that
\begin{equation*}
  \hat\chi_4=8 q_4 q_5^2-16 q_4^2 q_5-4q_6^2-4 \chi_1q_4 q_5 +\hat\chi_2q_5 -16 \chi_3 q_4,
\end{equation*}
where $\hat\chi_2:=\chi_2-\frac{1}{4} \chi_1^2 + 4\chi_3$ and $\hat\chi_4:=\chi_4 +2 \chi_1 \chi_3$, which can be
used instead of $\chi_2$ and $\chi_4$ as basic $\Ad$-invariant polynomials, restricted to $N$.  Using
(\ref{eq:atp_matrix_d4}), expressed in terms of the coordinates $\chi_1,\hat\chi_2,\chi_3,q_4,q_5$ and $q_6$ we find
that  the matrix $\Omega$ is indeed of the form (\ref{eq:omega}), with $c'=-1/8$, since
\begin{equation*}
\renewcommand{\arraystretch}{2}
  \begin{array}{rcl}
    \pb{q_4,q_5}&=& q_6 =\ds-\frac18\pp {\hat\chi_4}{q_6}, \\
    \pb{q_4, q_6}&=& 2q_4q_5-2q_4^2-\frac12 \chi_1q_4 +\frac18 \hat\chi_2=\ds\frac18\pp {\hat\chi_4}{q_5}, \\
    \pb{q_5,q_6}&=&-q_5^2+4 q_4 q_5+\frac12 \chi_1q_5  + 2 \chi_3 =\ds-\frac18 \pp {\hat\chi_4}{q_4}.
\end{array}
\end{equation*}
It follows easily from these formulas that the Milnor algebra is given by
$$\M(8 q_4 q_5^2-16 q_4^2q_5-4q_6^2)=\C[q_4,q_5,q_6]/(q_6,q_4(q_5-q_4),q_5(q_5-4q_4)),$$
so that $1$ and the coefficients $q_4,q_5$ and $q_4q_5$ of $\hat\chi_4$ form indeed a vector space basis for it.

\smallskip

\noindent
3) We finally consider the subregular orbit $\O_{sr}$ in $\Liesl_4$. This example is from \cite{DA}. This case was
also examined in \cite{SA} where it was shown that the slice, originally due to Arnold~\cite{AR}, belongs to the
set ${\mathcal N}_h$. It is the orbit of the nilpotent element
\begin{equation*}
  e=
  \begin{pmatrix}
    0 &1 &0 &0 \cr
    0&0& 1 & 0 \cr
    0&0&0& 0 \cr
    0&0&0&0
  \end{pmatrix}.
\end{equation*}
The transverse slice in Arnold's coordinates consists of matrices of the form
\begin{equation*}
  Q=
  \begin{pmatrix}
    0& 1& 0& 0 \cr
    0& 0 & 1 & 0 \cr
    q_1 &q_2 & q_3 & q_4 \cr
    q_5 & 0 &0 & -q_3
  \end{pmatrix}.
\end{equation*}
The basic Casimirs of the ATP-structure, as computed from the characteristic polynomial of $Q$, are
\begin{eqnarray*}
  \chi_1&=&q_2+q_3^2,\\
  \chi_2 &=& q_1+q_2 q_3, \\
  \chi_3&=&q_1 q_3 +q_4 q_5.
\end{eqnarray*}
If we solve the first two equations for the variables $q_1,\, q_2$ in terms of $\chi_1,\,\chi_2$ and
$q_3,\,q_4,\,q_5$, and substitute the result in $\chi_3$, then we find that
$$
  \chi_3=q_3^4+q_4 q_5-\chi_1q_3^2 +\chi_2 q_3.
$$
Using the explicit formulas for the ATP-structure, given in \cite{DA}, expressed in terms of the coordinates
$\chi_1,\chi_2,q_3,q_4$ and $q_5$, we find that the matrix $\Omega$ is indeed of the form (\ref{eq:omega}), with
$c'=1$, since
\begin{equation*}
\renewcommand{\arraystretch}{2}
  \begin{array}{rcl}
    \pb{q_3,q_4}&=& q_4 =\ds\pp {\chi_3}{q_5}, \\
    \pb{q_3, q_5}&=&-q_5=-\ds\pp {\chi_3}{q_4}, \\
    \pb{q_4,q_5}&=&4q_3^3-2\chi_1q_3+\chi_2=\ds\pp {\chi_3}{q_3}.
\end{array}
\end{equation*}
It can be read off from these formulas that the Milnor algebra is given by
$$\M(q_3^4+q_4q_5)=\C[q_3,q_4,q_5]/(q_4,q_5,q_3^3),$$
so that the coefficients $1,q_3$ and $q_3^2$ of $\chi_3$ form indeed a vector space basis for it.

\end{document}